\documentclass[a4paper, 12pt]{article}
\usepackage{enumerate, theorem}
\usepackage{amsmath, amsfonts, amssymb}
\usepackage[height=24cm, width=15cm]{geometry}
\usepackage[all]{xy}
\usepackage{mathrsfs, yfonts}

\newcommand{\parag}[1]{\paragraph{\sc{#1.}} }
\def\A{{\tilde{\mathcal{A}}}}

\newtheorem{thm}{Theorem}[subsection] 
\newtheorem{defn}[thm]{Definition}
\newtheorem{cor}[thm]{Corollary}
\newtheorem{prop}[thm]{Proposition}
\newtheorem{lemma}[thm]{Lemma}

\begin{document}

\title{Finite determination of regular (a,b)-modules.}

\author{Daniel Barlet\footnote{Barlet Daniel, Institut Elie Cartan UMR 7502  \newline
Nancy-Universit\'e, CNRS, INRIA  et  Institut Universitaire de France, \newline
BP 239 - F - 54506 Vandoeuvre-l\`es-Nancy Cedex.France.\newline
e-mail : barlet@iecn.u-nancy.fr}.}

\date{seconde version du 20 /08/07}

\maketitle

\markright{Finite determination of regular (a,b)-modules}

\tableofcontents

\parag{Summary}
The concept of (a,b)-module comes from the study the Gauss-Manin lattices of an isolated singularity of a germ of an holomorphic function. It is a very simple ''abstract algebraic structure'', but very rich, whose prototype is the formal completion of the Brieskorn-module of an isolated singularity.\\
The aim of this article is to prove a very basic theorem on regular (a,b)-modules showing that a given regular (a,b)-module is completely characterized by some ''finite order jet'' of its structure. Moreover a very simple bound for such a sufficient order is given in term of the rank and of two very simple invariants : the regularity order which count the number of times you need to apply \ $b^{-1}.a \simeq \partial_z.z$ \ in order to reach a simple pole (a,b)-module. The second invariant is the ''width'' which corresponds, in the simple pole case, to the maximal integral difference between to eigenvalues of \ $b^{-1}.a$ \ (the logarithm of the monodromy). \\
In the computation of examples this theorem is quite helpfull because it tells you at which power of \ $b$ \ in the expansions you may stop without loosing any information.

\newpage

\section*{Introduction.}

The concept of (a,b)-module comes from the study the Gauss-Manin lattices of an isolated singularity of a germ of an holomorphic function. It is a very simple ''abstract algebraic structure'', but very rich, whose prototype is the formal completion of the Brieskorn-module of an isolated singularity.\\
It appears that this structure induces an interesting approach in the study of singular points of linear differential systems (in one variable). As it will be apparent in this article, this point of view leads to study some finite type left  modules over the non-commutative \ $\mathbb{C}-$algebra generated by two variables \ $a,b$ :
$$\A : = \big\{ \sum_{\nu = 0}^{+\infty} \quad b^{\nu}.P_{\nu}(a) \quad  \big\} $$
where \ $P_{\nu} $ \ are in \ $ \mathbb{C}[z]$ \  and with the commutation relation \ $a.b - b.a = b^2$, assuming the continuity of left and right multiplication by \ $a$ \ for the \ $b-$adic topology of \ $\A$. Of course this commutation relation is satisfied by the ''standard model'' 
$$ a : = \times z, \quad b : = \int_0^z .$$
The aim of this article is to prove a very basic theorem on regular (a,b)-modules showing that a given regular (a,b)-module is completely characterized by some ''finite order jet'' of its structure. Moreover a very simple bound for such a sufficient order is given in term of the rank and of two very simple invariants : the regularity order which count the number of times you need to apply \ $b^{-1}.a \simeq \partial_z.z$ \ in order to reach a simple pole (a,b)-module. The second invariant is the ''width'' which corresponds in the simple pole case to the maximal integral difference between to eigenvalues of \ $b^{-1}.a$ \ (the logarithm of the monodromy). \\
In the computation of examples this theorem is quite helpfull because it tells you at which power of \ $b$ \ in the expansions you may stop without loosing any information.

\section{Basic properties.}

\subsection{Definition and examples.}

First recall the definition of an (a,b)-module.

\begin{defn}\label{(a,b)_module}
An {\bf (a,b)-module \ $E$} \ is a free finite type \ $\mathbb{C}[[b]]-$module with a \ $\mathbb{C}-$linear endomorphism \ $a : E \to E$ \ which is continuous for the \ $b-$adic topology of \ $E$ \ and satisfies \ $a.b - b.a = b^2$.\\
The {\bf rank of \ $E$}, denote by \ $rank(E)$,  will be the rank of \ $E$ \ as a \ $\mathbb{C}[[b]]-$module.
\end{defn}

\parag{Remarks}
\begin{enumerate}
\item Let \ $(e_1, \cdots, e_k)$ \ a \ $\mathbb{C}[[b]]-$basis of a free finite type  $\mathbb{C}[[b]]-$module. Then choosing arbitrarily elements \ $(\varepsilon_1, \cdots, \varepsilon_k)$ \ and defining \ $a.e_j = \varepsilon_j \quad \forall j \in [1,k]$ \ gives an (a,b)-module: the commutation relation implies \ $ a.b^n = b^n.a + n.b^{n+1} \quad \forall n \in \mathbb{N}$ \ so \ $a$ \ is defined on \ $\sum_{j =1}^k \quad \mathbb{C}[b].e_j$. The continuity assumption gives its (unique) extension.
\item There is a natural (a,b)-module associated to every algebraic linear differential system (see [B.95] p.42)
$$ Q(z).\frac{dF}{dz} = M(z).F(z), \quad Q \in \mathbb{C}[z], \quad M \in End(\mathbb{C}^n)\otimes_{\mathbb{C}} \mathbb{C}[z] .$$
\end{enumerate}

In the sequel of this article we shall mainly consider regular (a,b)-modules (see definition recalled below).  To try to convince the reader that the ''general'' (a,b)-module structure is interesting, let me quote the following result, which is quite elementary in the regular case, but which is not so easy in general.

\begin{thm}([B.95] th.1bis p.31)
Let \ $E$ \ be an (a,b)-module. Then the kernel and cokernel of ''a''  are finite dimensional.
\end{thm}

This result implies a general finiteness theorem for extensions of (a,b)-modules (see [B.95] and also section 1.3).

\bigskip

\begin{defn}\label{simple pole}
We shall say that an (a,b)-module \ $E$ \ has a {\bf simple pole} when the inclusion \ $a.E \subset b.E$ \ is satisfied.
\end{defn}

This terminology comes from the terminology of meromorphic connexions (see for instance [D.70]). 

\parag{Example}
 For any \ $\lambda \in \mathbb{C}$ \ define the simple pole  rank 1 \ (a,b)-module \ $E_{\lambda}$ \ as \ $E : = \mathbb{C}[[b]].e_{\lambda}$ \ where ''$a$'' is defined by the relation \ $a.e_{\lambda} = \lambda.b.e_{\lambda}$. $\hfill \square$
 
 \bigskip
 
 As an introduction to our main theorem, the reader may solve the following exercice by direct computation.
 \parag{Exercice}
For any \ $S \in \mathbb{C}[[b]]$ \ show that the simple pole (a,b)-module defined by  \ $E : = \mathbb{C}[[b]].e_S $ \ and \ $a.e_S = b.S(b).e_S$ \ is isomorphic to \ $E_{\lambda}$ \ with \ $\lambda = S(0)$ \\
(hint: begin by looking for \ $ \alpha_1 \in \mathbb{C}$ \ such that \ $(a - S(0).b)(e + \alpha_1.b.e) \in b^3.E$). $\hfill \square$

\bigskip

For a simple pole (a,b)-module, the linear map \ $b^{-1}.a : E \to E $ \ is well defined and induces an endomorphism \ $ f : = b^{-1}.a : E/b.E \to E/b.E$. For any \ $\lambda \in \mathbb{C}$ \ we shall denote by \ $\lambda_{min}$ \ the smallest eigenvalue of \ $f$ \ which is in \ $\lambda + \mathbb{Z}$. Then for \ $\lambda = \lambda_{min}- k$ \ with \ $k \in \mathbb{N}^*$ \ the bijectivity of the map \ $f - \lambda$ \ on \ $E/b.E$ \ implies easily its bijectivity on \ $E$ \ (see the exercice above). It gives then the equality 
$$ (a - \lambda.b).E = b.E.$$
Using this remark, it is not difficult to prove the following result from [B.93] (prop.1.3. p.11) that we shall use  later on.

\begin{prop}\label{sub min}
Let \ $E$ \ be a simple pole (a,b)-module, and let \ $\lambda \in \mathbb{C}$ \ and \ $\kappa \in \mathbb{N}$ \ such that \ $\lambda - \kappa \leq \lambda_{min} $. If \ $y \in E$ \ satisfies \ $(a - \lambda.b).y \in b^{\kappa+2}.E $ \ then there exists an unique \ $\tilde{y} \in E$ \ such that \ $(a - \lambda.b).\tilde{y} = 0 $ \ and \ $\tilde{y} - y \in b^{\kappa+1}.E$.
\end{prop}

An easy consequence of this proposition is that for an eigenvalue \ $\lambda$ \ of \ $f$ \ such that \ $\lambda = \lambda_{min}$ \ there always exists a non zero \ $x \in E$ \ such that \ $(a- \lambda.b).x = 0$. This gives an embedding of \ $E_{\lambda}$ \ in \ $E$. Remark also that if \ $E$ \ is a non zero simple pole (a,b)-module, such a \ $\lambda$ \ always exists. This leads to a rather precise description a of ''general'' simple pole (a,b)-module (see [B.93] th. 1.1 p.15). 

\bigskip

\begin{defn}\label{regular}
An (a.b)-module \ $E$ \ is  {\bf regular} \ when its saturation by \ $b^{-1}.a$ \ in \ $E[b^{-1}]$ \ is finitely generated on \ $\mathbb{C}[[b]]$.
\end{defn}

\bigskip

We shall denote \ $E^{\sharp}$ \ this saturation. It is a simple pole (a,b)-module and it is the smallest simple pole (a,b)-module containing \ $E$ \ in the sense that for any  (a,b)-linear morphism  \ $j : E \to F$ \ where \ $F$ \ is a simple pole (a,b)-module, there exists a unique (a,b)-linear extension \ $j^{\sharp} : E^{\sharp} \to F$ \ of \ $j$. \\

\smallskip

It is easy  to show that a regular (a,b)-module of rank 1  is isomorphic to some \ $E_{\lambda}$ \ for some \ $\lambda \in \mathbb{C}$. The classification of rank 2 regular (a,b)-module is not so obvious. We recall it here for a later use

\begin{prop}(see [B.93] prop.2.4 p. 34)\label{class}
The list of rank 2 regular (a,b)-modules is, up to isomorphism, the following :
\begin{enumerate}
\item \ $E_{\lambda} \oplus E_{\mu} $ \ for \ $(\lambda,\mu) \in \mathbb{C}^2/\frak{S}_2 $.
\item For any \ $\lambda \in \mathbb{C}$ \ and any \ $n \in \mathbb{N}$ \ let \ $E_{\lambda}(n)$ \ be the simple pole (a,b)-module with basis \ $(x,y)$ \ such that
$$  a.x = (\lambda + n).b.x + b^{n+1}.y \quad {\rm and} \quad  a.y = \lambda.b.y .$$
\item For any \ $(\lambda,\mu) \in \mathbb{C}^2/\frak{S}_2 $ \ let \ $E_{\lambda,\mu}$ \ the rank 2 regular (a,b)-module with basis \ $(y,t)$ \ such that
$$ a.y = \mu.b.y \quad {\rm and } \quad a.t = y + (\lambda-1).b.t .$$
\item For any \ $\lambda \in \mathbb{C}$, any \ $n \in \mathbb{N}^*$ \ and any \ $\alpha \in \mathbb{C}^*$ \ let \ $E_{\lambda, \lambda-n}(\alpha)$ \ be the rank 2 regular (a,b)-module with basis \ $(y,t)$ \ such that
$$ a.y = (\lambda-n).b.y \quad {\rm and} \quad a.t = y + (\lambda-1)b.t + \alpha.b^n.y $$
\end{enumerate}
\end{prop}

Note that the first two cases are simple pole (a,b)-modules.\\
 The saturation by \ $b^{-1}.a$ \  in case 3 is generated by \ $b^{-1}.y$ \ and \ $t$ \ as a \ $\mathbb{C}[[b]]-$module. It is isomorphic to \ $E_{\lambda-1} \oplus E_{\mu-1}$ \ for \ $\lambda \not= \mu$ \ and to \ $E_{\lambda-1}(0)$ \ for \ $\lambda = \mu$.\\
 The saturation  by \ $b^{-1}.a$ \  in case 4 is generated by \ $b^{-1}.y$ \ and \ $t$ \ as a \ $\mathbb{C}[[b]]-$module. It is isomorphic to \ $E_{\lambda-n-1}(n)$ \ for any non zero value of \ $\alpha$.
 
 \bigskip

To conclude this first section, let me recall also the theorem of existence of Jordan-H{\"o}lder sequences for regular (a,b)-module, which will be usefull in the induction in the proof of  our main  result .

\begin{thm}(see [B.93] th. 2.1 p.30)\label{J-H}
For any regular  rank k (a,b)-module \ $E$ \ there exists a sequence of sub-(a,b)-modules
$$ 0 = E^0 \subset E^1 \subset \cdots \subset E^{k-1} \subset E^k = E $$
such that for any \ $j \in [1,k]$ \ the quotient \ $E^j/E^{j-1}$ \ is isomorphic to \ $E_{\lambda_j}$. Moreover we may choose for \ $E^1$ \ any normal\footnote{normal means \ $E^1 \cap b.E = b.E^1$, so that \ $E/E^1$ \ is again free on \ $\mathbb{C}[[b]]$.}  rank 1 sub-(a,b)-module of \ $E$.\\
The number \ $\alpha(E) : = \sum_{j=1}^k \ \lambda_j $ \ is independant of the choice of the Jordan-H{\"o}lder sequence. It is given by the following formula
$$ \alpha(E) = trace\big(b^{-1}.a : E^{\sharp}/b.E^{\sharp} \to E^{\sharp}/b.E^{\sharp} \big) + \dim_{\mathbb{C}}(E^{\sharp}/E) .$$
\end{thm}

\subsection{The regularity order.}

\begin{defn}\label{order of regularity}
Let \ $E$ \ be a regular (a,b)-module. We define the {\bf regularity order of \ $E$} as the smallest integer \ $k \in \mathbb{N}$ \ such that the inclusion
\begin{equation*}
 a^{k+1}.E \subset \sum_{j=0}^k \quad a^j.b^{k-j+1}.E \tag{reg.}
 \end{equation*}
is valid. We shall note this integer \ $or(E)$.\\
We define also {\bf the index \ $\delta(E)$ \ of \ $E$} \ as the smallest integer \ $m \in \mathbb{N}$ \ such that \ $E^{\sharp} \subset b^{-m}.E $.
\end{defn}

\parag{Remarks}
\begin{enumerate}[i)]
\item The (a,b)-module \ $E$ \ has a simple pole if an only iff \ $or(E) = 0 $.
\item The inclusion \ (reg.) \ implies that \ $ (b^{-1}.a)^{k+1}.E \subset \Phi_k(E) : = \sum_{j=0}^k \quad (b^{-1}.a)^j.E$ \ and this implies that \ $\Phi_k(E)$ \ is stable by \ $b^{-1}.a$. So \ $\Phi_k(E)$ \ is a simple pole (a,b)-module contained in \ $b^{-k}.E \subset E[b^{-1}]$. This implies clearly the regularity of \ $E$. \\
For \ $k = or(E)$ \ we have \ $E^{\sharp} = \Phi_k(E) \subset b^{-k}.E$. So we have \ $\delta(E) \leq or(E)$.
\item As the quotient \ $b^{-k}.E/E$ \ is a finite dimensional \ $\mathbb{C}-$vector space, the quotient \ $E^{\sharp}/E$ \ is always a finite dimensional \ $\mathbb{C}-$vector space. $\hfill \square$
\end{enumerate}

The remark iii) shows that for a regular (a,b)-module \ $E$ \  there always exists a simple pole sub-(a,b)-module of \ $E$ \ which is  a finite codimensional vector space in \ $E$. This comes from the fact that for \ $k =  \delta(E)$ \ we have \ $ b^k.E^{\sharp} \subset E $ \ and that  \ $b^k.E^{\sharp}$ \ has a simple pole.

\parag{Example} The inequality \ $\delta(E) \leq or(E)$ \ may be strict for \ $or(E) \geq 2$. For instance the (a,b)-module of rank 3  with \ $\mathbb{C}[[b]]-$basis \ $e_1, e_2, e_3$ \ with \\
 $a.e_1 = e_2, \quad a.e_2 = b.e_3, \quad a.e_3 = 0 $ \ has index 1  and regularity order 2 : an easy computation gives that  a \ $\mathbb{C}[[b]]-$basis for \ $E^{\sharp} $ \ is given by \ $e_1, b^{-1}.e_2, b^{-1}.e_3$, and that a  \ $\mathbb{C}[[b]]-$basis for \ $E + b^{-1}.a.E$ \ is given by \ $e_1, b^{-1}.e_2, e_3$. $\hfill \square$

\begin{defn}\label{$E^b$}
Let \ $E$ \ be a regular (a,b)-module. The {\bf biggest simple pole sub-(a,b)-module of \ $E$} \ exists and has finite \ $\mathbb{C}-$codimension in \ $E$. We shall note it \ $E^b$.
\end{defn}

In general, for \ $k = \delta(E)$ \  the inclusion \ $b^k.E^{\sharp} \subset E^b $ \ is strict. For instance this is the case for \ $E_{\lambda,\mu} \oplus E_{\nu}$.

\bigskip

\begin{lemma}\label{inclusion} 
Let \ $E$ \ be a regular (a,b)-module. The smallest integer \ $m$ \ such we have \ $b^m.E \subset E^b$ \ is equal to \ $\delta(E)$.
\end{lemma}

\parag{Proof} Let \ $k : = \delta(E)$. Then \ $b^k.E^{\sharp}$ \ is a simple pole sub-(a,b)-module of \ $E$. So we have \ $b^k.E \subset b^k.E^{\sharp} \subset E^b$. Conversely, an inclusion \ $b^m.E \subset E^b $ \ gives \ $E \subset b^{-m}.E^b $. As \ $b^{-m}.E^b$ \ has a simple pole this implies \ $E^{\sharp} \subset b^{-m}.E^b \subset b^{-m}.E $. So \ $\delta(E) \leq m$.  $\hfill \blacksquare$

\parag{Examples} In the case 3 of the proposition \ref{class} \ $E^b$ \ is generated as a \ $\mathbb{C}[[b]]-$module by \ $y$ \ and \ $b.t$, so \ $E^b = b.E^{\sharp}$.\\
In case 4 we have also \ $E^b = b.E^{\sharp}$.

\bigskip

\begin{lemma}\label{ordre reg.}
Let \ $E$ \ be a regular (a,b)-module. For any exact sequence of (a,b)-modules
\begin{equation*}
 0 \to E' \to E \overset{\pi}{ \to} E''  \to 0 \tag{*}
\end{equation*}
we have  \ $or(E'') \leq or(E) \leq rank(E') +  or(E'') $.\\
As a consequence, the order of regularity of \ $E$ \ is at most \ $rank(E) - 1$ \ for any regular non zero (a,b)-module. 
\end{lemma}

\parag{Proof}  The inequality \ $or(E'') \leq or(E)$ \ is trivial because an inequality
$$  a^{k+1}.E   \subset \sum_{j=0}^{k} \ a^j.b^{k-j+1}.E $$
implies the same for \ $E''$ \ and, by definition, the best such integer \ $k$ \ is the order of regularity.\\
The crucial case is when \ $E'$ \ is of rank 1 . So we may assume that \ $E' \simeq E_{\lambda}$ \ for some \ $\lambda \in \mathbb{C}$ \ (see \ref{J-H} or  [B.93] prop.2.2 p.23). Let \ $k = or(E'')$. Then the inclusion 
\begin{equation*}
 a^{k+1}.E'' \subset \sum_{j=0}^{k} \ a^j.b^{k-j+1}.E'' \tag{1}
 \end{equation*}
 implies that 
 \begin{equation*}
 a^{k+1}.E \subset \sum_{j=0}^{k} \ a^j.b^{k-j+1}.E + b^l.E_{\lambda}  \tag{2}
 \end{equation*}
  for some \ $l \in \mathbb{N}$. In fact we can take for \ $l$ \  the smallest integer such that the generator \ $e_{\lambda}$ \ of \ $E_{\lambda}$ \ (defined up to \ $\mathbb{C}^*$ \ by the relation \ $a.e_{\lambda} = \lambda.b.e_{\lambda}$) satisfies \ $ b^l.e_{\lambda} \in \Psi_k  = \sum_{j=0}^{k} \ a^j.b^{k-j+1}.E $. \\ 
  Remark that this integer \ $l \geq 0$ \ is well defined because \ $b^{k+1}.e_{\lambda} \in \Psi_k$. Moreover, as \ $\Psi_k$ \ is a \ $\mathbb{C}[[b]]-$submodule of \ $E$, \ $b^l.e_{\lambda} \in \Psi_k$  \ implies \ $b^l.E_{\lambda} \subset \Psi_k$.\\
 Now, thanks to \ $(2)$ \  we have 
 \begin{equation*}
 a^{k+2}.E \subset \sum_{j=0}^k \ a^{j+1}.b^{k+1-j}.E \quad + a.b^l.E_{\lambda} \tag{3}
 \end{equation*}
 which gives
 \begin{equation*}
 a^{k+2}.E \subset \sum_{j=0}^{k+1} \ a^j.b^{k-j+2}.E  \tag{4}
 \end{equation*}
 because \ $a.b^l.E_{\lambda} = b.b^l.E_{\lambda} \subset b.\Psi_k $.\\
 This proves that \ $or(E) $ \ is at most \ $k + 1 = or(E'') + rank(E') $.\\
 Assume now that our inequality is proved for \ $E'$ \ of rank \ $p - 1$ \ and consider an exact sequence \ $(^*)$ \ with  \ $rank(E')$ \ equal \ $p \geq 2$. Let \ $E_{\lambda} \subset E'$ \ be a normal rank 1  sub-(a,b)-module of \ $E'$ \ (see \ref{J-H} or  [B.93]  prop.2.2 p.23 for a proof of the existence of such sub-(a,b)-module) and consider the exact sequence of (a,b)-modules (using the fact that \ $E_{\lambda}$ \ is also normal in \ $E$; see lemma 2.5 of [B.93])
 \begin{equation*}
 0 \to E'/E_{\lambda} \to E/E_{\lambda} \to E'' \to 0 
 \end{equation*}
 Using the induction hypothesis and the rank 1 case we get
 $$ or(E) \leq   or(E/E_{\lambda}) + 1 \leq p-1 + or(E'') + 1 = p + or(E'') .$$
 Now using an easy induction  (or a Jordan-H{\"o}lder sequence for \ $E$)  we obtain \ $or(E) \leq rank(E) - 1$ \ for any regular \ $E$. $\hfill \blacksquare$
 
 \parag{Remark} In the situation of the previous lemma we have \ $\delta(E') \leq \delta(E)$. This is a consequence of the obvious inclusion \ $ (E')^{\sharp} \subset E'[b^{-1}] \cap E^{\sharp}$ : assume that \ $x \in E'[b^{-1}] \cap E^{\sharp}$ ; then, for \ $k : = \delta(E)$ \ we have  \ $ b^k.x\in E'[b^{-1}] \cap E$ \ so that \ $b^{N+k}.x \in E'$ \ for \ $N$ \ large enough. As \ $E/E'$ \ has no $b-$torsion, we conclude that \ $b^k.x \in E'$. So our initial inclusion implies \ $\delta(E') \leq k$. $\hfill \square$
 
 \bigskip
 
 \subsection{Duality.}
 
 In this section we consider the associative and unitary \ $\mathbb{C}-$algebra
 $$ \tilde{\mathcal{A}} : = \big\{ \sum_0^{\infty} \  b^n.P_n(a) \quad {\rm with} \quad P_n \in \mathbb{C}[z] \big\} $$
with the commutation relation \ $a.b - b.a = b^2 $, and such that the left and right multiplications by \ $a$ \ are continuous for the \ $b-$adic topology\footnote{remark that for each \ $k \in \mathbb{N}$ \ $b^k.\A = \A.b^k$.} of \ $\A$.

\parag{The right structure as a commuting left-structure on \ $\A$} \quad \\
There exits an unique \ $\mathbb{C}-$linear (bijective) map \ $\theta : \A \to \A$ \ with the following properties
\begin{enumerate}[i)]
\item \ $\theta(1) = 1, \quad \ \theta(a) = a, \quad \theta(b) = - b $;
\item \ $\theta(x.y) = \theta(y).\theta(x) \quad \forall x,y \in \A $.
\item \ $\theta$ \ is continuous for the \ $b-$adic topology of \ $\A$
\end{enumerate}
The uniqueness is an easy consequence of iii) and  the fact that the conditions i) and ii) implies \ $\theta(b^p.a^q) = (-1)^p.a^q.b^p \quad \forall p,q \in \mathbb{N}$. Existence is then clear from the explicit formula deduced from this remark.\\
We  define a new structure of left \ $\A-$module on \ $\A$, {\bf called the \ $\theta-$structure} and denote by \ $x_*\square$,  by the formula
$$ x_*y = y.\theta(x) .$$
It is easy to see that this new left-structure on \ $\A$ \ commutes with the ordinary one and that with this \ $\theta-$structure \ $\A$ \ is still free of rank one as a left \ $\A-$module.

\begin{defn}\label{Hom}
Let \ $E$ \ be a (left) \ $\A-$module. On the \ $\mathbb{C}-$vector space \ $Hom_{\A}(E,\A)$ \ we define a  left  \ $\A-$module structure using the \ $\theta-$structure on \ $\A$. Explicitely this means that for \ $\varphi \in Hom_{\A}(E,\A)$ \ and \ $x \in \A$ \ we let
$$\forall e \in E \quad  (x.\varphi)(e) : = x_*\varphi(e) = \varphi(e).\theta(x).$$
We obtain in this way a left  \ $\A-$module that we shall still denote \ $Hom_{\A}(E,\A)$.
\end{defn}

It is clear that \ $E \to Hom_{\A}(E,\A)$ \ is a contravariant functor which is left exact in the category of left \ $\A-$modules. As every finite type left\ $\A-$module has a resolution of length \ $\leq 2$ \ by free finite type modules ( see [B.95] cor.2 p.29), we shall denote by \ $Ext^i_{\A}(E, \A), i \in [0,2]$ \ the right derived functors of this functor. They are finite type left \ $\A-$modules when \ $E$ \ is finitely generated because \ $\A$ \ is left noetherian (see [B.95] prop.2 p.26).

\bigskip

Any (a,b)-module is a left \ $\A-$module. They are characterized by the existence of special simple resolutions.

 \begin{lemma}\label{Resol.}
 Let \ $M$ \ be a \  $(p,p)$ \ matrix with entries in the ring \ $\mathbb{C}[[b]] $. Then the left \ $\tilde{\mathcal{A}}-$linear map \ $ Id_p.a - M : \tilde{\mathcal{A}}^p \to \tilde{\mathcal{A}}^p$ \ given by \\ $$^tX : = (x_1, \cdots,x_p) \  \to  \  ^tX.(Id_p.a - M)$$
 is injective. Its cokernel is the (a,b)-module \ $E$ \ given as follows : \\ 
 $E$ \ has a \ $\mathbb{C}[[b]]$ \ base \ $e : = (e_1, \cdots, e_p)$ \ and \ $a$ \ is defined by the two conditions
 \begin{enumerate}
 \item \ $a.e : = M(b).e $ ;
 \item the left action of \ $a$ \ is continuous for the \ $b-$adic topology of \ $E$.
 \end{enumerate}
 Any (a,b)-module is obtained in this way and so, as a \ $\A-$left-module, has a resolution of the form
 \begin{equation*}
  0 \to \A^p \  \overset{^t\square.(Id_p.a - M)}{\longrightarrow} \ \A^p \to E \to 0 . \tag{@}
  \end{equation*}
  \end{lemma}
  
  \parag{Proof} First remark that for \ $x \in \A$ \ the condition \ $x.a \in b.\A $ \ implies \ $x \in b.\A $. Now let us prove, by induction on \ $n \geq 1 $, that, for any \ $(p,p)$ \ matrix \ $M$ \ with entries in \ $\mathbb{C}[[b]]$ \ the condition \ $ ^tX.(Id_p.a - M) = 0 $ \ implies \ $^tX \in b^n.\A^p$.\\
   For \ $n =1$ \ this comes from the previous remark. Let assume that the assertion is proved for \ $n \geq 1$ \ and consider an \ $X \in \A^p$ \ such that \ $ ^tX.(Id_p.a - M) = 0 $. Using the induction hypothesis we can find \ $Y \in \A^p$ \ such that \ $X = b^n.Y$. Now we obtain, using \ $a.b^n = b^n.a + n.b^{n+1}$ \  and the fact that \ $\A$ \ has no zero divisor, the relation
   $$ ^tY(Id_p.a -(M + n.Id_p.b)) = 0 $$
   and using again our initial remark we conclude that \ $Y \in b.\A^p $ \ so \ $X \in b^{n+1}.\A^p$.\\
   So such an \ $X$ \ is in \ $\cap_{n\geq 1} \ b^n.\A^p = (0)$. \\
   The other assertions of the lemma are obvious. $\hfill \blacksquare$
   
   \bigskip
   
   We recall now a construction given in [B.95] which allows to compute more easily the vector spaces \ $Ext^i_{\A}(E,F)$ \ when \ $E,F$ \ are (a,b)-modules 
   
   \begin{defn}\label{Hom(a,b)}
   Let \ $E,F$ \ two (a,b)-modules. Then the \ $\mathbb{C}[[b]]-$module \ $Hom_b(E,F)$ \ is again a free and finitely generated \ $\mathbb{C}[[b]]-$module. Define on it an (a,b)-module structure in the following way.
   \begin{enumerate}
   \item First change the sign of the action of \ $b$. So \ $S(b) \in \mathbb{C}[[b]]$ \ acts as \ $\check{S}(b) = S(-b)$.
   \item Define \ $a$ \ using the linear map
   $ \Lambda : Hom_b(E,F) \to Hom_b(E,F)$ \ given by \ $\Lambda(\varphi)(e) = \varphi(a.e) - a.\varphi(e)$.
   \end{enumerate}
   We shall denote \ $Hom_{a,b}(E,F)$ \ the corresponding (a,b)-module.
   \end{defn}
   The verification that \ $\Lambda(\varphi)$ \ is \ $\mathbb{C}[[b]]-$linear and that \ $\Lambda.\check{b} - \check{b}.\Lambda = \check{b}^2$ \ are easy (and may be found in [B.95] p.31).
   
  \parag{Remark} In {\it loc. cit.} we defined the (a,b)-module structure on \ $Hom_{a,b}(E,F)$ \ with opposite  signs for \ $a$ \ and \ $b$. The present convention is better because it fits with the usual definition of the formal adjoint of a differential operator : \ $z^* = z$ \ and \ $(\partial/\partial z)^* = -\partial/\partial z$. $\hfill \square$
  
  \bigskip

   The following lemma is also proved in {\it loc.cit.}
   
   \begin{lemma}\label{hom}
    Let \ $E,F$ \ two (a,b)-modules. Then there is a functorial isomorphism of \ $\mathbb{C}-$vector spaces 
     $$ H^i\Big( Hom_{a,b}(E,F) \overset{a}{\to}Hom_{a,b}(E,F) \Big) \to Ext^i_{\A}(E, F)\quad \forall i \geq 0 .$$
     Here the map \ $a$ \ of the complex \ $Hom_{a,b}(E,F) \overset{a}{\to}Hom_{a,b}(E,F) )$ \ is equal to the \ $\Lambda$ \ defined above which is, by definition,  the operator \ $''a''$ \ of the (a,b)-module \ $Hom_{a,b}(E,F)$.
     \end{lemma}
     
     Now the following corollary of the lemma  \ref{Resol.} gives that the two natural ways of defining the dual of an (a,b)-module give the same answer.  
   
     \begin{cor}\label{Dual Reg.}
   Let \ $E$ \ an (a,b)-module. There is a functorial isomorphism of (a,b)-modules between the following two (a,b)-modules constructed as follows : 
   \begin{enumerate}
   \item \ $Ext^1_{\A}(E,\A)$ \ with the \ $\A-$structure defined by the \ $\theta-$structure of \ $\A$.
   \item \ ${\it Hom}_{a,b}(E, E_0) $ \ where \ $E_0 : = \A \big/ \A.a$.
   \end{enumerate}
   \end{cor}
   
   \parag{Proof} Using a free resolution \ $(@)$ \ of \ $E$ \ deduced from a \ $\mathbb{C}[[b]]-$basis \\
    $e : = (e_1, \cdots, e_p)$ \ we obtain, by the previous lemma, an exact sequence
    \begin{equation*}
  0 \to \A^p \  \overset{(Id_p.a - ^tM).\square}{\longrightarrow} \ \A^p \to Ext^1_{\A}(E, \A) \to 0 . \tag{@@}
  \end{equation*}
  of left \ $\A-$modules where \ $\A^p$ \ is endowed with its \ $\theta-$structure. Writing the same exact sequence with the ordinary left-module structure of \ $\A^p$ \ gives
   \begin{equation*}
  0 \to \A^p \  \overset{^t\square.(Id_p.a - ^t\check{M})}{\longrightarrow} \ \A^p \to Ext^1_{\A}(E, \A) \to 0 . \tag{@@ bis}
  \end{equation*}
  where \ $^t\check{M}(b) : = \  ^tM(-b) $.\\
  Denote by \ $e^* : = (e_1^*, \cdots, e_p^*)$ \ the dual basis of \ $Hom_{\mathbb{C}[[b]]}(E,E_0)$. By definition of the action of \ $a$ \ on \ ${\it Hom}_{a,b}(E, E_0) $ \ we get, if \ $\omega$ \ is the class of 1 in \ $E_0$ :
  $$ (a.e_i^*)(e_j) = e_i^*(a.e_j) - a.e_i^*(e_j) = e_i^*(\sum_{h=1}^p \ m_{j,h}.e_h) - a.\delta_{i.j}.\omega = \check{m}_{j,i}.\omega $$ 
  because \ $a.\omega = 0$ \ in \ $E_0$, and the definition of the action of \ $b$ \ on \ ${\it Hom}_{a,b}(E, E_0) $. So we have \ $ a.e^* = ^t\check{M}.e^*$ \ concluding the proof. $\hfill \blacksquare$
  
  \begin{defn}\label{Dual} For any (a,b)-module \ $E$ \ the {\bf dual of \ $E$}, denoted by \ $E^*$, is the (a,b)-module \ $Ext^1_{\A}(E,\A) \simeq {\it Hom}_{a,b}(E, E_0) $.
  \end{defn}
  
  Of course, for any \ $\A-$linear map \ $f : E \to F$ \ between two (a,b)-modules we have an\ $\A-$linear ''dual'' map \ $f^* : F^* \to E^* $.\\
  It is an easy consequence of our previous description of \ $Ext^1_{\A}(E,\A)$ \ that we have a functorial isomorphism \ $(E^*)^*\to E$.
 
    \parag{Examples}
  \begin{enumerate}
  \item For each \ $\lambda \in \mathbb{C}$ \ we have \ $(E_{\lambda})^* \simeq E_{-\lambda}$.
  \item For \ $(\lambda,\mu) \in \mathbb{C}^2$ \ we have \ $E_{\lambda,\mu}^* \simeq E_{-\mu+1,-\lambda+1}$.
  \item Let \ $E$ \ be  the rank two simple pole (a,b)-module \ $E_1(0)$ \ defined by \ $a.e_1 = b.e_1 + b.e_2$ \ and \ $a.e_2 = b.e_2$.  Then its dual is isomorphic to \ $E_{-1}(0)$.\\
     It is also an elementary exercice to show the following isomorphisms :
    $$E_1(0) \simeq \mathbb{C}[[z]] \oplus  \mathbb{C}[[z]].Log z \quad {\rm and} \quad E_{-1}(0) \simeq  \mathbb{C}[[z]]\frac{1}{z^2} \oplus  \mathbb{C}[[z]].\frac{Log z}{z^2}$$ 
     with \ $a : = \times z$ \ and \ $ b : = \int_0^z $.
  \end{enumerate}
  
  \begin{prop}\label{Dualite et regularite}
  For any exact sequence of (a,b)-modules
  $$ 0 \to E' \overset{u}{\to} E \overset{v}{\to} E'' \to 0 $$
  we have an exact sequence of (a,b)-modules
  $$ 0 \to (E'')^* \overset{v^*}{\to}  E^* \overset{u^*}{\to}  (E')^* \to 0 .$$
  If \ $E$ \ is a simple p\^ole (a,b)-module, $E^*$ \ has a simple pole. \\
  For any regular (a,b)-module \ $E$ \ its dual \ $E^*$ \ is regular. Moreover, if \ $E^b$ \ and \ $E^{\sharp}$ \ are respectively the biggest simple pole submodule of \ $E$ \ and the saturation of \ $E$ \ by \ $b^{-1}.a$ \ in \ $E[b^{-1}]$, we have
  $$ (E^{\sharp})^* \simeq (E^*)^b \quad {\rm and} \quad  (E^b)^* \simeq (E^*)^{\sharp}.$$
  \end{prop} 
   
  \parag{Proof}
  The first assertion is a direct consequence of the vanishing of \ $Ext^i_{\A}(E, \A)$ \ for \ $i = 0, 2$, for any (a,b)-module and the long exact sequence for the ''Ext''.\\
  The condition that \ $E$ \ has a simple pole is equivalent to the fact that for any choosen basis \ $e$ \ of \ $E$ \ the matrix \ $M$ \ has its coefficients in \ $b.\A = \A.b$. Then this remains true for \ $^t\check{M}$.\\
  To prove the regularity of  \ $E^*$ \ when \ $E$ \ is  regular, we shall use induction on the rank of \ $E$. The rank 1 case is obvious because we have a simple pole in this case. Assume that the assertion is true for \ $ rank < p $ \ and consider a \ $rank = p$ \   regular (a,b)-module \ $E$. Using the theorem \ref{J-H} we have an exact sequence of (a,b)-modules
  $$ 0 \to E_{\lambda} \to E \to F \to 0  $$
  where \ $F$ \ is regular of rank \ $p-1$. This gives a short exact sequence
      $$ 0 \to F^* \to E^* \to E_{-\lambda} \to 0 $$
  and the regularity of \ $F^*$ \ and of \ $E_{-\lambda}$ \ implies the regularity of \ $E^*$.\\
  Now the inclusions \ $ E^b \subset E \subset E^{\sharp}$ \ gives exact sequences
  \begin{align*}
  & 0 \to Ext^1_{\A}(E/E^b, \A) \to E^* \to (E^b)^* \to Ext^2_{\A}(E/E^b, \A) \to 0 \\
  &  0 \to Ext^1_{\A}(E^{\sharp}/E, \A)  \to (E^{\sharp})^* \to E^* \to  Ext^2_{\A}(E^{\sharp}/E, \A) \to 0
  \end{align*}
  and the next lemma will show that the \ $Ext^1_{\A}(V, \A) = 0 $ \ for any \ $\A-$module which is a finite dimensional vector space, and also the finiteness (as a vector space) of \ $Ext^2_{\A}(V, \A)$.
  This implies  that we have, for any regular (a,b)-module, the inclusions
  $$ E^* \subset  (E^b)^*\quad {\rm and} \quad (E^{\sharp})^* \subset E^* .$$
  They imply, thanks to the fact that \ $(E^b)^*$ \ and \ $(E^{\sharp})^*$ \ have simple poles,
  $$  (E^*)^{\sharp} \subset (E^b)^* \quad {\rm and} \quad (E^{\sharp})^* \subset (E^*)^b.$$
    But the inclusion \ $(E^*)^b \subset E^* $ \ gives
     $$ E = (E^*)^*\subset ((E^*)^b)^* \subset ((E^{\sharp})^*)^* = E^{\sharp}$$
     and the minimality of \ $E^{\sharp}$ \ gives \ $((E^*)^b)^* = E^{\sharp}$ \ because \ $((E^*)^b)^*$ \ has a simple pole and contains \ $E$. Dualizing again  gives \ $(E^{\sharp})^* \simeq (E^*)^b$. 
     The last equality is obtained in a similar way from \ $ E^* \subset (E^*)^{\sharp}$. $\hfill \blacksquare$

  \begin{lemma}\label{Finite dim. modules}
  Let \ $V$ \ be an \ $\A-$module of finite dimension over \ $\mathbb{C}$. Then we have \ $Ext^i_{\A}(V,\A) = 0 $ \ for \ $i = 0,1$ \ and \ $Ext^2_{\A}(V,\A)$ \ is again a \ $\A-$module (via the \ $\theta-$structure of \ $\A$) which is a finite dimensional vector space. Moreover it has the same dimension than \ $V$ \ and there is a canonical \ $\A-$module isomorphism
  $$ Ext^2_{\A}(Ext^2_{\A}(V,\A),\A) \simeq V .$$
  \end{lemma}
  
  \parag{proof} We begin by proving the first assertion of the lemma for the special case \ $ V_{\lambda} : = \A \big/ \A.(a-\lambda) + \A.b$ \ for any  \ $\lambda \in \mathbb{C}$. Let us show that we have the free resolution
  $$ 0 \to \A \overset{\alpha}{\to} \A^2 \overset{\beta}{\to} \A \to V_{\lambda} \to 0 $$
  where \ $\alpha(x) : = (x.b, -x.(a-b-\lambda)), \quad \beta(u,v) : = u.(a-\lambda) + v.b $. The map \ $\alpha$ \ is clearly injective and \ $\beta(\alpha(x)) = x.(b.a - \lambda.b - (a-b-\lambda).b) = 0$. If we have \ $\beta(u,v) = 0$ \ then \ $u \in \A.b $; let \ $u = x.b$. Then we get
  $$ x.(a-b-\lambda).b + v.b = 0 \quad {\rm and \ so} \quad v = -x.(a-b-\lambda) .$$
  This gives the exactness of our resolution.\\
  Now the \ $Ext^i_{\A}(V_{\lambda},\A)$ \ are given by the cohomology of the complex
  $$ 0 \to \A   \overset{\beta^*}{\to}  \A^2  \overset{\alpha^*}{\to} \A \to 0 .$$
  The map \ $\beta^*(x) = ((a-\lambda).x, b.x)$ \ and \ $\alpha^*(u,v) = b.u - (a-b-\lambda).v $ \ are \ $\A-$linear for the \ $\theta-$structure of \ $\A$. Clearly \ $\beta^*$ \ is injective and \ $\alpha^*( \beta^*(x)) \equiv 0$. If \ $\alpha^*(u,v) = 0 $ \ set \ $v = b.y$ \ and conclude that \ $u = (a-\lambda).y$. This gives the vanishing of the \ $Ext^i$ \ for \ $i = 0,1$. The \ $Ext^2$ \ is the cokernel of \ $\beta^*$ \ which is easily seen to be isomorphic to \ $V_{\lambda}$.\\
  Consider now any finite dimensional \ $\A-$module \ $V$ \ over \ $\mathbb{C}$. We make an induction on \ $\dim_{\mathbb{C}}(V)$ \ to prove the vanishing of the \ $Ext^i$ \ for \ $i = 0,1$ \ and the assertion on the dimension of the \ $Ext^2$.\\
   The \ $\dim V = 1 $ \ case is clear because reduced to the case \ $V = V_{\lambda}$ \ for some \ $\lambda \in \mathbb{C}$. Assume that the case \ $\dim V = p$ \ is proved, for \ $p \geq 1$ \ and consider some \ $V$ \ with \ $\dim V = p +1$. Then \ $Ker \, b$ \ is not \ $\{0\}$ \ and is stable by \ $a$. Let \ $\lambda \in \mathbb{C}$ \ an eigenvalue of \ $a$ \ acting on \ $ Ker \, b$. Then a eigenvector generates in \ $V$ \ a sub-$\A-$module isomorphic to \ $V_{\lambda}$. \\
  The exact sequence of \ $\A-$modules
  $$ 0 \to V_{\lambda} \to V \to W \to 0 $$
  where \ $W : = V\big/ V_{\lambda}$ \ has dimension \ $p$ \ allows us to conclude, looking at the long exact sequence of  Ext .\\
  The last assertion follows from the remark that we produce a free resolution of \ $Ext^2_{\A}(V,\A)$ \ by taking \ $Hom_{\A}(-,\A)$ \ of a free (length two, see [B.97]) resolution of \ $V$ \ because of the already proved vanishing of the \ $Ext^i$ \  for \ $i = 0,1$. Doing this again gives back the initial resolution (remark that we use here that the \ $\theta\circ\theta-$structure on \ $Hom_{\A}(Hom_{\A}(\A,\A),\A)$ \ is the usual left structure on \ $\A$). $\hfill \blacksquare$
  
  \begin{cor}\label{Sym. Spec.}
 For a simple pole (a,b) module  \ $E$ \ denote by \ $S(E)$ \ the spectrum of \ $b^{-1}.a$ \ acting on \ $E/b.E$. Then we have
 $$ S(E^*) = - S(E).$$
 \end{cor}
 
 \parag{Proof} We make an induction on the rank of \ $E$. In rank $1$ \ the result is clear because we have \ $E \simeq E_{\lambda}$ \ for some \ $\lambda \in \mathbb{C}$, and \ $S(E_{\lambda}) = \{\lambda \}$. But we know that \ $E_{\lambda}^* = E_{-\lambda}$.\\
 Assume the assertion proved for any rank \ $p \geq 1$ \ simple pole (a,b)-module, and consider \ $E$ \ with rank \ $p +1$.  Using theorem \ref{J-H}, there exists \ $\lambda \in \mathbb{C}$ \ and an exact sequence  (a,b)-modules
 $$ 0 \to E_{\lambda} \to E \to F \to 0 $$
 where \ $rank(F)= p$ \ and where \ $F$ \ has a simple pole (because a quotient of a simple pole (a,b)-module has a simple pole !). The exact sequence of vector spaces 
 $$  0 \to E_{\lambda}/b.E_{\lambda}  \to E/b.E \to F/b.F \to 0 $$
 shows that \ $S(E) = S(F) \cup \{\lambda \}$. Now proposition \ref{Dualite et regularite} gives the exact sequence
 $$ 0 \to F^* \to E^* \to  E_{-\lambda} \to 0 $$
 which implies, as before, $S(E^*) = S(F^*) \cup \{-\lambda\} $. The induction hypothesis \ $S(F^*) = - S(F)$ \ allows to conclude. $\hfill \blacksquare$
 
 \bigskip
  
  \begin{lemma}\label{Duality for extension}
  For any pair of (a,b)-modules\ $E$ \ and\ $F$ \ there is a canonical isomorphism of vector spaces
  $$ D : Ext^1_{\A}(E,F) \to Ext^1_{\A}(F^*, E^*)  $$
 associated to the correspondance between 1-extensions (i.e. short exact sequences)
 $$ (0 \to F \to G \to E \to 0) \overset{D}{\to} (0 \to E^* \to G^* \to F^* \to 0).$$
 \end{lemma}
 
 \parag{Proof} We have a obvious isomorphism of \ $\mathbb{C}[[b]]-$modules\footnote{but be carefull with the \ $ b \to \check{b}$ !}
 $$ I :  Hom_b(E,F) \to Hom_b(Hom_b(F,E_0), Hom_b(E,E_0))\simeq Hom_b(F^*,E^*) $$
 because  \ $E_0 \simeq \mathbb{C}[[b]]$ \ as  a\ $\mathbb{C}[[b]]-$module. But recall that \ $Ext^1_{\A}(E,F)$ \ (resp. \ $Ext^1_{\A}(F^*,E^*)$) \ is the cokernel of the \ $\mathbb{C}-$linear map ''$a$'' defined on \ $ Hom_b(E,F) $ \ by the formula
 $$ (a.\varphi )(x) = \varphi(a.x) - a.\varphi(x) $$
 So it is enough to check that the isomorphism \ $I$ \ commutes with ''$a$'' in order to get an isomorphism between the cokernels of ''$a$'' in these two spaces.\\
 Let \ $\varphi \in Hom_b(E,F)$ and \ $\xi \in F^*$. Then \ $I(\varphi)(\xi) = \varphi\circ \xi $. 
 So, for \ $x \in  E$ \ we have (using \ $\Lambda$ \ to avoid too many ''$a$'')
 \begin{align*}
 & \Lambda(I(\varphi)(\xi) = I(\varphi)(a.\xi) - a.(I(\varphi)(\xi)) \\
 &  \Lambda(I(\varphi)(\xi)(x) = (\varphi \circ \xi )(a.x) - a.\xi(\varphi(x)) - \big(\xi(\varphi(a.x)) - a.\xi(\varphi(x))\big) \\
 & \quad\quad \quad\quad    = \big[(\Lambda(\varphi))\circ \xi \big](x) = I(\Lambda(\varphi))(x).
 \end{align*}
 So \ $\Lambda\circ I = I \circ \Lambda$. The map \ $I$ \ gives an isomorphism of complexes\\
 $$\xymatrix{Hom_{a,b}(E,F) \ar[r]^{\Lambda} \ar[d]^I & Hom_{a,b}(E,F)  \ar[d]^I \\
  Hom_{a,b}(F^*,E^*) \ar[r]^{\Lambda} &Hom_{a,b}(F^*,E^*) }$$
  and this conclude the proof, using lemma \ref{hom}. $\hfill \blacksquare$
  
  \bigskip
  
  For an (a,b)-module \ $E$ \ and an integer \ $m \in \mathbb{N}$ \ it is clear that \ $b^m.E$ \ is again an (a,b)-module. This can be generalize for any \ $m \in \mathbb{C}$.
  
  \begin{defn}
  For any (a,b)-module \ $E$ \ and any complex number \ $m \in \mathbb{C}$ \ define the (a,b)-module \ $b^m.E$ \ as follows :  as an \ $\mathbb{C}[[b]]-$module we let \ $b^m.E \simeq E \simeq \mathbb{C}[[b]]^{rank(E)}$;  the operator \ $a$ \ is defined as \ $a + m.b$.
  \end{defn}

  Precisely, this means that if \ $(e_1, \cdots, e_k)$ \ is a  \ $\mathbb{C}[[b]]-$basis of \ $E$ \ such that we have \ $a.e = M(b).e $ \ where \ $M \in End(\mathbb{C}^p) \otimes_{\mathbb{C}} \mathbb{C}[[b]]$, the (a,b)-module \ $b^m.E$ \ admit a basis, denote by \ $(b^m.e_1, \cdots, b^m.e_k)$, such that the operator \ $a$ \ is defined by the relation \ $a.(b^m.e) : = (M(b) + m.b.Id_k).(b^m.e) $.\\
  Remark that for \ $m \in \mathbb{N}$ \ this notation is compatible with the preexisting one, because of the relation \ $a.b^m = b^m.(a + m.b)$.\\
  For any \ $m \in \mathbb{N}$ \ there exists a canonical (a,b)-morphism
  $$ b^m.E \to  E  $$
  which is an isomorphism of \ $b^m.E$ \ on \ $Im(b^m : E \to E)$. But remark that the map \ $b^m : E \to E$ \ is not $a-$linear (but the image is stable by \ $a$).\\
  For any \ $m \in \mathbb{N}$ \ there is also a canonical (a,b)-morphism
  $$ E \to b^{-m}.E $$
  which induces an isomorphism of \ $E$ \ on \ $Im(b^m : b^{-m}.E \to b^{-m}.E)$. So we may write, via this canonical identification,  $b^m.(b^{-m}.E) = E$.\\
  It is easy to see that for any \ $m, m' \in \mathbb{C}$ \ we have a natural isomorphism
  $$ b^{m'}.(b^m.E) \simeq b^{m+m'}.E \quad {\rm and \ also} \quad b^0.E \simeq E .$$
  
  \parag{Remark} It is easy to show that for any \ $m \in \mathbb{C}$ \ there exists an unique $\mathbb{C}-$algebra automorphism 
   $$ \eta_m : \A \to \A \quad {\rm such \  that } \quad \eta(1) = 1, \eta(b) = b \quad {\rm and} \quad \eta(a) = a +m.b .$$
   Using this automorphism, one can define a left \ $\A-$module  \ $b^m.F$ \ for any left \ $\A-$module \ $F$ \ and any \ $m \in \mathbb{C}$. This is, of course compatible with our definition in the context of (a,b)-modules. $\hfill \square$
   
   \bigskip
   
   The behaviour of the correspondance \ $E \to b^m.E$ \  by duality is given by the following easy lemma; the proof is left as an exercice.
   
   \begin{lemma}
   For any (a,b)-module \ $E$ \ and any \ $m \in \mathbb{C}$ \ there is natural (a,b)-isomorphism
   $$ (b^m.E)^*\to b^{-m}.E^* .$$
   \end{lemma}
   
   \bigskip
   
   The following corollary of the lemma  \ref{inclusion} and the proposition \ref{Dualite et regularite} allows to show that duality preserves the index.
   
   \begin{lemma}\label{or dual}
   Let \ $E$ \ be a regular (a,b)-module. Then  we have \ $\delta(E^*) = \delta(E)$.
   \end{lemma}
   
   \parag{Proof} By definition \ $\delta(E)$ \ is the smallest integer \ $k \in \mathbb{N}$ \  such that \ $E^{\sharp} \subset b^{-k}.E$.\\
 Now \ $E^{\sharp} \subset b^{-m}.E$ \  implies by duality that \ $b^m.E^* \subset (E^*)^b$. So,  by lemma \ref{inclusion}, we have \ $m \geq \delta(E^*)$. This proves that \ $\delta(E) \leq \delta(E^*)$ \ and we obtain the equality by symetry. $\hfill \blacksquare$
 
\parag{Remark} Duality does not preserve the order of regularity : in the example given before the definition \ref{$E^b$} we have \  $or(E) = 2$ \ and \ $or(E^*) = 1$. $\hfill \square$

\bigskip
   
   Let us conclude this section by an easy exercice.
   
   \parag{Exercice}
   For any (a,b)-modules  \ $E, F$ \ and any \ $\lambda \in \mathbb{C}$ \ there are natural (a,b)-isomorphisms
   \begin{enumerate}
   \item  \ $b^{\lambda}.E_{\mu} \simeq E_{\lambda+\mu}$.
    \item \ $b^{\lambda}.Hom_{a,b}(E, F) \simeq  Hom_{a,b}(b^{-\lambda}.E, F) \simeq Hom_{a,b}(E, b^{\lambda}.F)$. 
    \item Then deduce from the previous isomorphisms that \ $Hom_{a,b}(E, E_{\lambda}) \simeq b^{-\lambda}.E^*$, and \ $Ext^1_{\A}(E, E_{\lambda}) \simeq E^*/(a + \lambda.b).E^*$.
    \end{enumerate}

\subsection{Width of a regular (a,b)-module.}

\parag{Notation} For a complex number \ $\lambda$ \ we shall note by \ $\tilde{\lambda}$ \ is class in \ $\mathbb{C}\big/\mathbb{Z}$. We shall order elements in each class \ $modulo \  \mathbb{Z}$ \ by its natural order on real parts. $\hfill \square$

\begin{defn}\label{largeur}
Let \ $E$ \ be a regular (a,b)-module and let \ $\tilde{\lambda} \in \mathbb{C}\big/\mathbb{Z}$. We define the following complex numbers :
\begin{align*}
& \tilde{\lambda}_{min}(E) : = \inf \{\lambda \in \tilde{\lambda} /\exists \ {\rm a \ non \ zero \ morphism} \quad E_{\lambda} \rightarrow E  \} \\
& \tilde{\lambda}_{max}(E) = \sup \{ \lambda \in \tilde{\lambda} / \exists  \ {\rm a \ non \ zero \ morphism} \quad  E \to E_{\lambda}  \} \\
& L_{\tilde{\lambda}}(E) = \tilde{\lambda}_{max}(E) - \tilde{\lambda}_{min}(E) \in \mathbb{Z} \\
& L(E) = \sup \{ \tilde{\lambda} \in  \mathbb{C}/\mathbb{Z} \  /  \ L_{\tilde{\lambda}}(E) \}
\end{align*}
with the following conventions :
 \begin{align*}
 & \inf \{\emptyset \} = +\infty, \ \sup \{\emptyset \} = - \infty \quad {\rm and} \\
 & -\infty - \lambda = - \infty \quad \forall \lambda \in ]-\infty, +\infty] \\
 & +\infty - \lambda = + \infty \quad \forall \lambda \in [-\infty, +\infty[. 
  \end{align*}
We shall call \ $L(E)$ \  {\bf the width of\ $E$}.
\end{defn}

\parag{Remarks}
\begin{enumerate}
\item A non zero morphism \ $E_{\lambda} \to E $ \ is necessarily injective. Either its image is a normal submodule in \ $E$ \ or there exists an integer \ $k \geq 1$ \ and a morphism \ $E_{\lambda-k} \to E$ \ whose image is normal an contains the image of the previous one.
\item In a dual way, a non zero morphism \ $E \to E_{\lambda}$ \  has an image equal to \ $b^k.E_{\lambda} \simeq E_{\lambda+k}$, where \ $k \in \mathbb{N}$.
\item A non zero morphism \ $E_{\lambda} \to E_{\mu}$ \ implies that \ $\lambda $ \ lies in \ $\mu + \mathbb{N}$. It is possible that for some \ $E$ \ we have \ $\tilde{\lambda}_{max}(E) < \tilde{\lambda}_{min}(E) $. For instance this is the case for the rank 2 regular (a,b)-module \ $E_{\lambda,\mu}$ \ from \ref{class}. So the width of a regular but not simple pole (a,b)-module  is not necessarily a non negative integer. 
\item Let \ $E$ \ and \ $F$ \ be regular (a,b)-modules. If there is a surjective morphism \ $ E \to F$ \ then for all \ $\tilde{\lambda} \in \mathbb{C}\big/\mathbb{Z}$ \ we have \ $\tilde{\lambda}_{max}(E) \geq \tilde{\lambda}_{max}(F)$. \\
If there is an injective morphism \ $E' \to E$ \ then for all \ $\tilde{\lambda} \in \mathbb{C}\big/\mathbb{Z}$ \ we have \ $\tilde{\lambda}_{min}(E) \leq \tilde{\lambda}_{min}(E')$.

\item Every submodule of \ $E$ \ isomorphic to \ $E_{\lambda}$ \ is contained in \ $E^b$. So we have \ $\tilde{\lambda}_{min}(E) = \tilde{\lambda}_{min}(E^b)$, for every regular (a,b)-module \ $E$ \ and every  \ $\tilde{\lambda} \in \mathbb{C}\big/\mathbb{Z}$.
\item In a dual way, every morphism \ $E \to E_{\lambda}$ \ extends uniquely to a morphism \ $E^{\sharp} \to  E_{\lambda}$ \ with the same image. So for every regular (a,b)-module \ $E$ \ and every  \ $\tilde{\lambda} \in \mathbb{C}\big/\mathbb{Z}$, we get \ $\tilde{\lambda}_{max}(E) = \tilde{\lambda}_{max}(E^{\sharp})$. $\hfill \square$
\end{enumerate}

\begin{lemma}\label{obvious}
\begin{enumerate}
\item Let \ $E$ \  a simple pole (a,b)-module and let \ $S(E)$ \ denotes the spectrum of the linear map \  $b^{-1}.a : E/b.E \to E/b.E$, we have
\begin{equation*}
 \tilde{\lambda}_{min}(E) = \inf \{ \lambda \in S(E) \cap \tilde{\lambda} \} \quad {\rm and} \quad   \tilde{\lambda}_{max}(E) = \sup \{ \lambda \in S(E) \cap \tilde{\lambda} \} \tag{@}
 \end{equation*}
\item For any regular (a,b)-module \ $E$ \ we have
$$ \widetilde{(-\lambda)}_{max}(E^*) = - \tilde{\lambda}_{min}(E) \quad\qquad \widetilde{(-\lambda)}_{min}(E^*) = -  \tilde{\lambda}_{max}(E).$$
This implies \ $ L_{-\tilde{\lambda}}(E^*) = L_{\tilde{\lambda}}(E)\  \forall \tilde{\lambda} \in \mathbb{C}/\mathbb{Z}$, and so \ $L(E^*) = L(E)$.
\item For any regular (a,b)-module \ $E$ \ and any \ $\tilde{\lambda} \in \mathbb{C}\big/\mathbb{Z}$ \ we have equivalence between
$$\tilde{\lambda}_{min}(E)  \not= +\infty \quad {\rm and} \quad \tilde{\lambda}_{max}(E) \not= -\infty .$$
\end{enumerate}
\end{lemma}

\parag{Proof} Let \ $E$ \ be a simple pole (a,b)-module. We have already seen (in proposition \ref{sub min})  that if \ $\lambda \in S(E)$ \ is minimal in its class modulo $1$, there exists a non zero \ $x \in E$ \ such that \ $a.x = \lambda.b.x$. This implies that \ $\tilde{\lambda}_{min} \leq \inf \{ \lambda \in S(E) \cap \tilde{\lambda} \}$. But the opposite inequality is obvious, so the first part of  (@) is proved.\\
Using corollary \ref{Sym. Spec.} and the result already obtained  for \ $E^*$ \ gives
$$ \widetilde{(-\lambda)}_{min}(E^*) = \inf \{-\lambda \in S(E^*) \cap \widetilde{(-\lambda)} \} = - \sup \{\lambda \in S(E) \cap \tilde{\lambda}\}.$$
So for \ $\mu = \sup \{\lambda \in S(E) \cap \tilde{\lambda}\}$ \ we have an exact sequence of (a,b)-modules
$$ 0 \to E_{-\mu} \to E^* \to F \to 0 $$
and by duality, a surjective map \ $E \to E_{\mu}$. This implies \ $\tilde{\lambda}_{max} \geq \mu$. As, again, the opposite inequality is obvious, the second part of (@) is proved.\\
Let us prove now the relations in 2. \\
Remark first that these equalities are true for a simple pole (a,b)-module because of \ $(@)$ \ and corollary \ref{Sym. Spec.}.\\
For any regular (a,b)-module \ $E$ \ we know that 
 $$\tilde{\lambda}_{min}(E) =  \tilde{\lambda}_{min}(E^b) =  \inf \{ \lambda \in S(E^b) \cap \tilde{\lambda} \} \quad {\rm and} \quad  \widetilde{(-\lambda)}_{max}(E^*) = \widetilde{(-\lambda)}_{max}((E^*)^{\sharp}).$$
  But we have
   $$  \widetilde{(-\lambda)}_{max}((E^*)^{\sharp}) = \sup\{-\lambda \in S((E^*)^{\sharp}) \cap \widetilde{(-\lambda)} \} = - \inf\{\lambda \in S((E^*)^{\sharp})^*\cap \tilde{\lambda} \} $$ 
   because \ $(E^*)^{\sharp}$ \  has a simple pole, using corollary \ref{Sym. Spec.}. So we obtain
   $$ \widetilde{(-\lambda)}_{max}(E^*) = - \tilde{\lambda}_{min}(E^b) = - \tilde{\lambda}_{min}(E)$$
   because \ $(E^*)^{\sharp})^* = E^b$ \ (see proposition \ref{Dualite et regularite}).\\  
The second relation is analoguous.\\
The equivalence in 3  is obvious in the simple pole case using \ $(@)$.\\
 The general case is an easy consequence using \ $E^b, E^{\sharp}$ : if \ $\tilde{\lambda}_{min}(E) \not= +\infty$ \ so is \ $\tilde{\lambda}_{min}(E^{\sharp})$ \ because \ $E \subset E^{\sharp}$. Then \ $\tilde{\lambda}_{max}(E^{\sharp})\not= -\infty$ \ and so is \ $\tilde{\lambda}_{max}(E)$. The converse is analoguous using \ $E^b$. \hfill $\blacksquare$

\parag{Remarks}
\begin{enumerate}
\item If \ $E$ \ has a simple pole, we have   \ $L_{\tilde{\lambda}}(E) \geq 0$ \ or \ $L_{\tilde{\lambda}}(E) = - \infty$ \ for any \ $\tilde{\lambda}$ \ in \ $\mathbb{C}/\mathbb{Z}$. So \ $L(E)$ \ is always \ $\geq 0$.
\item In cases 1 and 2 of the proposition \ref{class} the formula \ $(@)$ \ gives the values of \ $\tilde{\lambda}_{min}$ \ and \ $\tilde{\lambda}_{max}$ \ for any \ $\tilde{\lambda} \in \mathbb{C}/\mathbb{Z}$.\\
For the remaining cases we can compute these numbers using the fact that we already know the corresponding \ $E^b$ \ and \ $E^{\sharp}$ \ and the remark 5 and 6  before the preceeding lemma. $\hfill \square$
\end{enumerate}

\bigskip

\begin{prop}\label{Induction largeur}
Let \ $E$ \ be a regular (a,b)-module and let \ $\tilde{\lambda} \in \mathbb{C}\big/\mathbb{Z}$. Assume that \ $ \lambda = \tilde{\lambda}_{min}(E) < +\infty $. Consider an exact sequence of (a,b)-modules
$$ 0 \to E_{\lambda} \to E \overset{\pi}{\to} F \to 0 .$$
Then we have for all  \ $\tilde{\mu} \in \mathbb{C}/\mathbb{Z}$ \ the inequality
\begin{equation*}
 L_{\tilde{\mu}}(F) \leq  L_{\tilde{\mu}}(E) + 1. \tag{i}
 \end{equation*}
 \end{prop}
 
 \parag{Proof} As \ $\tilde{\mu}_{max}(F) \leq \tilde{\mu}_{max}(E)$ \ for any \ $\mu \in \mathbb{C}$ \ it is enough to prove that we have \ $ \tilde{\mu}_{min}(E) \leq \tilde{\mu}_{min}(F) +1$ \ for all \ $\tilde{\mu} \in \mathbb{C}/\mathbb{Z}$.
 
  Let begin by the case of \ $\tilde{\mu} = \tilde{\lambda}$. We want to show the inequality
 \begin{equation*}
 \tilde{\lambda}_{min}(F) \geq \lambda -1 \tag{ii}
 \end{equation*}
 Let  \ $E_{\lambda - d} \hookrightarrow F$ \ with \ $d \geq 0$. The rank 2 (a,b)-module  \ $G : = \pi^{-1}(E_{\lambda - d})$ \ is contained in \ $E$, so \ $\lambda = \tilde{\lambda}_{min}(G)$. We have
  the exact sequence of (a,b)-modules
 $$ 0 \to E_{\lambda} \to \pi^{-1}(E_{\lambda - d}) \overset{\pi}{\to} E_{\lambda-d} \to 0. $$
 Now let us  compare \ $G$ \ with the list in proposition \ref{class}.\\
  If \ $G$ \ is in case 1, we have \ $E_{\lambda-d} \subset G$ \ so \ $d = 0$ \ because \ $\lambda = \tilde{\lambda}_{min}(G)$.\\
 If \ $G$ \ is in case 2, we have \ $\lambda - d = \lambda + n$ \ with \ $n \in \mathbb{N}$, so \ $d = 0$.\\
 If \ $G$ \ is in case 3, we have  \ $G \simeq E_{\lambda, \lambda+k}$ \ with \ $k \in \mathbb{N}$. Then  the theorem \ref{J-H} gives \ $2\lambda - d = 2\lambda +k -1 $ \ and so \ $d = 1-k \leq 1$.\\
 If \ $G$ \ is in case 4, we have \ $G \simeq E_{\lambda,\lambda+n}(\alpha)$. Again theorem \ref{J-H} gives \ $2\lambda - d = 2\lambda+n-1 $ \ so \ $d = 1-n \leq 0$ \ because \ $n \in \mathbb{N}^*$. So \ $d = 0$.\\
  We conclude that we always have \ $d \leq 1$ \ and this proves (ii).\\
   
    \smallskip
   
For \ $\tilde{\mu} \not= \tilde{\lambda}$ \ let us prove now  the following inequality :
\begin{equation*}
\tilde{\mu}_{min}(F) \leq  \tilde{\mu}_{min}(E) \leq \tilde{\mu}_{min}(F) +1 . \tag{iii}
\end{equation*}
Consider an injective morphism  \ $E_{\mu} \to E$ \ with \ $\mu = \tilde{\mu}_{min}(E)$. The restriction of \ $\pi$ \ to \ $E_{\mu}$ \ is injective and so it gives \  $ \tilde{\mu}_{min}(E) \geq \tilde{\mu}_{min}(F)$. Assume now that we have an injective morphism  \ $E_{\mu'} \hookrightarrow F$ \ with  \ $\mu' = \tilde{\mu}_{min}(F)$, and consider the rank 2 (a,b)-module \ $\pi^{-1}(E_{\mu'})$. Using the proposition \ref{class} where only cases 1 or 3 are possible now, it can be easily check that (iii) is satisfied. $\hfill \blacksquare$

\parag{Remarks} 
\begin{enumerate}
\item In the situation of the previous proposition we have either \ $\tilde{\lambda}_{min}(E) \geq \tilde{\lambda}_{max}(E)$ \ or \ $\tilde{\lambda}_{max}(E)= \tilde{\lambda}_{max}(F)$ : Assume that we have  \ $\lambda < \lambda' : = \tilde{\lambda}_{max}(E)$. Then there exists a surjective morphism \ $ q : E \to E_{\lambda'}$, and, as the restriction of \ $q$ \ to  \ $E_{\lambda}$ \ is zero, the map \ $q$ \ can be factorized and gives a surjective morphism \ $\tilde{q} : F \to E_{\lambda'}$. So we get \ $\tilde{\lambda}_{max}(E)\leq \tilde{\lambda}_{max}(F)$, and the desired equality thanks to the preceeding lemma.
\item We shall use later that in the situation of the previous proposition we have the inequality \ 
$ \tilde{\lambda}_{max}(F) \leq \lambda + L(E) . \hfill  \square$
\end{enumerate}

\begin{cor}\label{maj. lambda(max)(F)}
In the situation of the previous proposition we have the inequality  \ $L(E) + rank(E) \geq L(F) + rank(F)$. So this integer is always positive for any non zero regular (a,b)-module.
\end{cor}

\parag{Proof} As  the rank  1 case is obious, an easy induction on the rank of \ $E$ \ using the propositions \ref{sub min}  and \ref{Induction largeur} gives the proof. $\hfill \blacksquare$

\parag{Examples}
\begin{enumerate}
\item The (a,b)-module
$$  J_k(\lambda) : = \A\big/\A.(a -(\lambda+k-1).b)(a -(\lambda+k-2).b)\cdots (a-\lambda.b)$$
which has rank \ $k$, satisfies \ $\lambda_{max} = \lambda$ \ and \ $\lambda_{min} = \lambda+k-1$. So its width is \ $L(J_k(\lambda)) = -k + 1 $.\\
To understand easily the (a,b)-module \ $J_k(\lambda)$ \ the reader may use the following alternative definition of it : there is a \ $\mathbb{C}[[b]]-$basis  \ $(e_1, \cdots, e_k)$ \ in which the action of ''$a$'' is given by 
$$ a.e_1 = e_2 + \lambda.b.e_1,\quad  a.e_2 = e_3 + (\lambda+1).b.e_2, \cdots, a.e_k =(\lambda+ k-1).b.e_k .$$
\item The rank 2 (a,b)-module \ $E_{\lambda} \oplus E_{\lambda+n}$ \ has width \ $n$. This shows that, despite the fact that the width is always bigger than \ $- rank(E) + 1$, the width may be arbitrarily big , even for a rank 2  regular (a,b)-module. $\hfill \square$
\end{enumerate}

 \section{Finite determination of regular (a,b)-modules.}

\subsection{Some more preliminaries.}

 \begin{lemma} \label{finite det. reg.}
 Let \ $E$ \ be a regular (a,b)-module of index \ $\delta(E) = k$. For \ $N \geq k+1$ \ the quotient map \ $ q_N : E \to E\big/b^N.E $ \ induces a bijection between simple pole sub-(a,b)-modules \ $F$ \ containing \ $b^k.E^{\sharp}$ \ and sub \ $\A-$modules \ $\mathcal{F} \subset  E\big/b^N.E $ \ satisfying the following two conditions
 \begin{enumerate}[i)]
 \item \ $ a.\mathcal{F} \subset b.\mathcal{F}$ ;
 \item  \ $b^k.E^{\sharp}\big/b^N.E \subset \mathcal{F}$.
 \end{enumerate}
 \end{lemma}
 
 \parag{Proof} It is clear that if \ $F$ \ is a simple pole sub-(a,b)-module of \ $E$ \ containing \ $b^k.E^{\sharp}$ \  the image \ $\mathcal{F} : = q_N(F)$ \ is a \ $\A-$submodule of \ $E\big/b^N.E $ \ such that i) and  ii) are fullfilled. Conversely, if a \ $\A-$submodule \ $\mathcal{F}$ \ satisfies i) and ii), let \ $F : = q_N^{-1}(\mathcal{F})$. Of course, \ $F$ \ is a sub-(a,b)-module of \ $E$ \ and contains \ $b^k.E^{\sharp}$. The only point to see is that \ $F$ \ has a simple pole. If \ $x \in F$ \ then \ $a.q_N(x) \in b.\mathcal{F}$ \ so \ $a.x \in b.F + b^N.E$. As \ $N \geq k+1$ \ we may write \ $ a.x = b.y + b.z $ \ with \ $y \in F $ \ and \ $z \in b^{N-1}.E \subset  b^k.E^{\sharp} \subset F $. This completes the proof. $\hfill \blacksquare$.
 
 \parag{Remarks} 
 \begin{enumerate}
 \item we may replace \ $b^k.E^{\sharp}$ \ by \ $b^k.E$ \ in the second condition imposed on \ $F$ \ and \ $\mathcal{F}$ : if a simple pole (a,b)-submodule \ $F$ \ contains \ $b^k.E$ \ it contains \ $b^k.E^{\sharp}$ \ by definition of \ $E^{\sharp}$. This allows to avoid the use of \ $E^{\sharp}$ \ in the previous lemma.
 \item The biggest \ $\mathcal{F}$ \ satisfying i) and ii) corresponds to \ $E^b$. So we may recover \ $E^b$ \ from the quotient \ $E\big/b^N.E$ \ for \ $N \geq \delta(E)+1$. $\hfill \square$
 \end{enumerate}
 
 \begin{cor}\label{det. finie ordre reg.}
  Let \ $E$ \ be a regular (a,b)-module of order of regularity \ $k$. Fix \ $N \geq k+1$ \ and assume that we has an isomorphism of \ $\A-$modules 
   $$ \varphi : E\big/b^N.E \to  E'\big/b^N.E'  $$
   where \ $E'$ \ is an (a,b)-module. Then \ $E'$ \ is regular and has order of regularity \ $k$. Moreover we  have the equality \ $\varphi(E^b\big/b^N.E) = (E')^b\big/b^N.E' $.
  \end{cor}
  
  \parag{Proof} As \ $k$ \ is the order of regularity of \ $E$ \ we have \ $a^{k+1}.E \subset \sum_{j=0}^{k} \   a^j.b^{k-j+1}E $. The inequality \ $N \geq k+1$ \ gives \ $a^{k+1}.E\big/b^N.E \subset \sum_{j=0}^{k} \  a^j.b^{k-j+1}E\big/b^N.E $, and this is also true for \ $E'\big/b^N.E'$, and implies \ $a^{k+1}.E' \subset \sum_{j=0}^{k} \  a^j.b^{k-j+1}E' $. So the order of regularity of \ $E'$ \ is at most \ $k$. We conclude that it is exactly \ $k$ \ by symetry.\\
  The last stament comes from the second remark above, as \ $or(E) \geq \delta(E)$. $\hfill \blacksquare$
  
\subsection{Finite determination for a rank one extension.}

 \begin{lemma}\label{existence N}
Let \ $E$ \ be an (a,b)-module et  fix a complex number \ $\lambda$. There exists \ $N(E,\lambda) \in \mathbb{N}$ \ such that for any \ $N \geq N(E,\lambda)$ we have the following inclusion :
  $$b^N.E \subset (a - \lambda.b).E .$$
\end{lemma}

\parag{Proof} With the $b-$adic topology, \ $E$ \ is a Frechet space. The \ $ \mathbb{C}-$linear map \ $a - \lambda.b : E \to E $ \ is continuous. The finiteness theorem of [B.95], theorem 1.bis p.31 gives that the kernel and cokernel of this map are finite dimensional vector spaces. So the subspace \ $(a - \lambda.b).E$ \ is closed in \ $E$. This statement corresponds to the equality
\begin{equation*}
 \cap_{N \geq 0} \big[(a - \lambda.b).E + b^N.E \big] = (a- \lambda.b).E  \tag{@}
 \end{equation*}
 But the images of the subspaces \ $b^N.E$ \ in the finite dimensional vector space \\
  $E \big/ (a- \lambda.b).E$ \ is a decreasing sequence. So it is stationnary, and, as the intersection is \ $\{0\}$ \ thanks to \ $(@)$, the result follows. $\hfill \blacksquare$

  \begin{prop}\label{4}
   Let \ $F$ \ be an (a,b)-module and \ $\lambda$ \ a complex number. Consider a short exact sequence of (a,b)-modules
  \begin{equation*}
   0 \to E_{\lambda} \overset{\alpha}{\to} E \overset{\beta}{\to} F \to 0  \tag{$@@$}
   \end{equation*}
  where \ $E_{\lambda} : = \A\big/\A.(a - \lambda.b) $. Then, for any \ $N \geq N(F^*, -\lambda)$, the extension \ $(@@)$ \ is uniquely determined by the following extension of \ $\A-$modules which are finite dimensional vectors spaces
  \begin{equation*}
   0 \to  E_{\lambda}\big/b^N.E_{\lambda}  \overset{\alpha}{\to} E\big/b^N.E \overset{\beta}{\to} F\big/b^N.F  \to 0  \tag{$@@_N$}   \end{equation*}
 obtained from \ $(@@)$ \ by ''quotient by \ $b^N$''.
\end{prop}

\parag{Comments} This statement needs some more explanations. Denote by \ $K_N$ \ the kernel of the obvious map (forget "a") 
 $$ob_N : Ext^1_{\A}(F/b^N.F, E_{\lambda}/b^N.E_{\lambda}) \to Ext^1_b(F/b^N.F, E_{\lambda}/b^N.E_{\lambda})$$
 where \ $Ext^1_b(-,-)$ \ is a short notation for \ $Ext^1_{\mathbb{C}[[b]]}(-,-)$. The  short exact sequence correspondance \ $(@@) \to (@@_N)$ \ gives a map 
 $$\delta_N : Ext^1_{\A}(F, E_{\lambda}) \to Ext^1_{\A}(F/b^N.F, E_{\lambda}/b^N.E_{\lambda}) $$
 whose  range  lies in \ $K_N$, because the \ $\mathbb{C}[[b]]-$exact sequence \ $(@@)$ \ is split as \ $F$ \ is \ $\mathbb{C}[[b]]-$free, and so is the exact sequence \ $(@@_N)$. The precise signification of the previous proposition is that for \ $N \geq N(F^*, -\lambda)$ \ the map \ $\delta_N$ \ is a \ $\mathbb{C}-$linear  isomorphism between the vector spaces  \ $Ext^1_{\A}(F, E_{\lambda})$ \ and \ $K_N$. $\hfill \square$
 
  \parag{Proof} As a first step to realize the map \ $\delta_N$ \ let us consider  the following commutative  diagramm of complex vector spaces, deduced from the exact sequences of \ $\A-$modules:
  \begin{align*}
 &  0 \to E_{\lambda + N} \to E_{\lambda} \to E_{\lambda}/b^N.E_{\lambda} \to 0  \\
  & 0 \to b^N.F \to F \to F/b^N.F \to 0 
  \end{align*}

 $$ \xymatrix{ Ext^1( F/ b^N.F,E_{\lambda+N}) \ar[d] \ar[r] &  Ext^1(F,E_{\lambda+N})  \ar[d]^{\alpha}  \ar[r] &  Ext^1(b^N.F,E_{\lambda +N})  \ar[d] \\
  Ext^1( F/ b^N.F, E_{\lambda}) \ar[d] \ar[r] & Ext^1(F, E_{\lambda})  \ar[d]^{\beta} \ar[r]^u  & Ext^1(b^N.F,E_{\lambda})  \ar[d]^v \\
  Ext^1(F/ b^N.F,E_{\lambda}/b^N.E_{\lambda}) \ar[r]^>>>>{\gamma} &  Ext^1(F, E_{\lambda}/b^N.E_{\lambda}) \ar[r]^w  & Ext^1(b^N.F, E_{\lambda}/b^N.E_{\lambda})} $$
  
  We have the following propreties : 
  
  \begin{enumerate}
 \item The surjectivity of the map \ $\beta$ \  is consequence of the vanishing of the vector space \ $Ext^2_{\A}(F, E_{\lambda +N})$ \ thanks to the proposition \ref{Dualite et regularite}.
  \item the vanishing of the composition \ $u\circ v $ \  is consequence of lemma  \ref{hom} and of the fact that the restriction map 
  $$ Hom_b(F, E_{\lambda}) \to Hom_b(b^N.F, E_{\lambda}) \to Hom_b(b^N.F, E_{\lambda}/b^N.E_{\lambda})$$ 
   is obviously zero.
   \item So the map \ $w$ \ is zero and \ $\gamma$ \ is surjective.
   \item The kernel of \ $\gamma$ \ is given by the image of the injective map
 $$ \partial : Hom_{\A}(b^N.F, E_{\lambda}/b^N.E_{\lambda}) \hookrightarrow Ext^1_{\A}(F/ b^N.F,E_{\lambda}/b^N.E_{\lambda}) .$$ 
 This  is  a consequence of the vanishing of the map
 $$ Ext^0_{\A}(F, E_{\lambda}/b^N.E_{\lambda}) \to Ext^0_{\A}(b^N.F, E_{\lambda}/b^N.E_{\lambda}).$$
 \end{enumerate}
 Let us show now that for \ $N \geq N(F^*, -\lambda)$ \ the map \ $\alpha$ \ is zero.  Using again the isomorphisms given by the lemma \ref{hom}, \ $\alpha$ \ is induced by the obvious map \ $ Hom_b(F, b^N.E_{\lambda}) \to Hom_b(F, E_{\lambda})$, whose image is \ $b^N.Hom_b(F, E_{\lambda})$. Denote respectively by \ $G$ \ and \ $H$ \  the (a,b)-modules given by \ $Hom_b(F, b^N.E_{\lambda})$ \ and \ $Hom_b(F, E_{\lambda})$ \ with the action of "$a$" defined by \ $\Lambda$ (see \ref{hom}). Then we have the following commutative diagramm
   $$ \xymatrix{ G \ar[r]^{i} \ar[d] & H \ar[d] \\ G/a.G \ar[r] \ar[d]^{\simeq} & H/a.H \ar[d]^{\simeq}\\
   Ext^1_{\A}(F, b^N.E_{\lambda}) \ar[r]^{\alpha} & Ext^1_{\A}(F, E_{\lambda})}$$
   and the image of \ $i$ \ is \ $b^N.H$. So the map \ $\alpha$ \ will be zero as soon as \ $b^N.H \subset a.H$ \ and this is fullfilled for \ $N \geq  N(H,0) = N(F^*, -\lambda) $. This last equality coming from the isomorphisms
   $$H/a.H \simeq Ext_{\A}^1(F, E_{\lambda}) \simeq Ext_{\A}^1( E_{-\lambda}, F^*) \simeq F^*/(a+\lambda.b).F^*$$
   see the exercice concluding \S 1.3.
   
   \smallskip
   
  Consider now the commutative diagramm
   
   $$ \xymatrix{\quad & 0 \ar[d] & K_N \ar[d]^i &Ext_{\A}^1(F, E_{\lambda})  \ar[d]^{\beta} \ar[l]_{\hat{\delta}_N} \ar[ld]_{\delta_N}  \\
    0 \ar[r] & Hom_{\A}(b^N.F, E_{\lambda}/b^N.E_{\lambda}) \ar[r]^{\partial}\ar[d]^{ob_N} &  Ext^1_{\A}(F/ b^N.F,E_{\lambda}/b^N.E_{\lambda})\ar[r]^>>>>>{\gamma} \ar[d]^{ob_N}& Ext^1_{\A}(F, E_{\lambda}/b^N.E_{\lambda})\\ \quad & Hom_b(b^N.F, E_{\lambda}/b^N.E_{\lambda})\ar[r]^{\simeq} &  Ext^1_b(F/ b^N.F,E_{\lambda}/b^N.E_{\lambda})& \quad }$$
    The surjectivity of \ $\beta$ \ implies that the map \ $i\circ\gamma $ \ is surjective ( we know that the extensions in the image of \ $\delta_N$ \ comes from \ $K_N$, so \ $\delta_N$ \ factors in \ $\hat{\delta}_N\circ i$).\\
    We have \ $i(K_N) \cap Im(\partial_N) = (0)$ \ because \ $ob_N$ \ is injective on \ $Im(\partial_N)$. \\
    So \ $i$ \ induces an isomorphism of vector spaces from \ $K_N$ \ to 
     $$ Ext^1_{\A}(F/ b^N.F,E_{\lambda}/b^N.E_{\lambda})/ Im(\partial_N)\overset{\gamma}{ \simeq} Ext^1_{\A}(F, E_{\lambda}/b^N.E_{\lambda}) \overset{\beta^{-1}}{\simeq} Ext_{\A}^1(F, E_{\lambda}).$$
     This completes the proof . $\hfill \blacksquare$
     
     \bigskip
     We shall need some bound for the integer \ $N(F^*, -\lambda)$ \ which appears in the previous proposition for the proof of our theorem.
      
      \begin{lemma}
      Let \ $G$ \ be a regular (a,b)-module and let \ $\mu \in \mathbb{C}$. A sufficient condition on \ $N \in \mathbb{N}$ \ in order to have the inclusion \ $ b^N.G \subset (a-\mu.b).G$ \ is  
      $$N \geq \mu -\tilde{\mu}_{min}(G) + \delta(G) + 2.$$
      \end{lemma}
      
      \parag{Proof} As we know that \ $ \tilde{\mu}_{min}(G^b) =  \tilde{\mu}_{min}(G)$ , for \ $M \in \mathbb{N}$, the assumption \ $M > \mu -  \tilde{\mu}_{min}(G)$ \  implies that \ $(a - (\mu-M).b).G^b = b.G^b$ \ (see the remark before proposition \ref{sub min}). By  definition of the index of \ $G$ \ we have  \ $b^{\delta(G)}.G \subset G^b $. Combining both gives 
      $$ b^{M+\delta(G)+1}.G \subset b^M.(a -(\mu-M).b).G = (a-\mu.b).b^M.G \subset (a-\mu.b).G.$$
      Now let \ $N = M + \delta(G) +1$ ; a sufficient condition on the integer \ $N$ \  is now \ $N \geq \mu -\tilde{\mu}_{min}(G) + \delta(G) + 2.$ \ $\hfill \blacksquare$

     \begin{cor}\label{3}
     A sufficient condition for \ $N \geq N(F^*,-\lambda)$ \ in the situation of prop. \ref{4} in the regular case is that \ $N \geq or(E) + L(E) + rank(E) + 1 $.
     \end{cor}
     
     Remark that the inequality \ $L(E) + rank(E) \geq 1$ \ for any non zero regular \ $E$ \ implies that we have  \ $or(E) + L(E) + rank(E) + 1 \geq or(E)  +  2$.
     
     \parag{Proof} We apply the previous lemma with \ $ F^*= G$ \ and \ $\mu = -\lambda = -\tilde{\lambda}_{min}(E)$. The conclusion comes now from the following facts :
     \begin{enumerate}
     \item \ $-\widetilde{(-\lambda)}_{min}(F^*) = \tilde{\lambda}_{max}(F) \leq \lambda + L(E) $ \ this last inequality is proved in \ref{Induction largeur}.
          \item \ $ \delta(F^*) = \delta(F) \leq or(F) \leq or(E) $ \ proved in \ref{or dual} and \ref{ordre reg.} $\hfill \blacksquare$
          \end{enumerate}

     \subsection{The theorem.}
     
     \begin{thm}\label{finite det. thm}
     Let \ $E$ \ be a regular (a,b)-module. There exists an integer \ $N(E)$ \ such that for any (a,b)-module \ $E'$,  any integer \ $ N \geq N(E)$ \ and any \ $\A-$isomorphism
     \begin{equation*}
      \varphi :  E/b^N.E \to  E'/b^N.E'  \tag{1}
      \end{equation*}
     there exists an unique \ $\A-$isomorphism \ $ \Phi :  E \to E' $ \ inducing the given \ $\varphi$.\\
     Moreover the choice \ $N(E)= N_0(E) : =  or(E) + L(E) + rank(E) + 1 $ \ is possible.
     \end{thm}
     
     \parag{Remarks}
     \begin{enumerate}
     \item It is easy to see that for a rank 1 regular (a,b)-module the integer 2  is the best possible.
     \item In our final lemma \ref{final} we show that the integer given in the theorem is optimal for the rank $k$ (a,b)-module \ $J_k(\lambda)$, (defined in the lemma), for any \ $k \in \mathbb{N}^*$.
     \item For the rank 2 (a,b)-modules \ $E_{\lambda, \mu}$ \ the integer  given by the theorem  is \ $or(E) + L(E) + 2 + 1  =  3$ \ is again optimal, as it can be shown in the same maner that in our final lemma. 
      \item  For the rank 2 simple pole (a,b)-module \ $E_{\lambda}(0)$ \ the integer given by the theorem is  \ $ 3 = L(E) + rank(E) + 1 $ \ and  the best possible is \ $2$ : the action of \ $b^{-1}.a$ \ on \ $E/b.E$ \ which is determined by \ $E/b^2.E$ \ characterizes this rank 2 regular (a,b)-module in the classification given in proposition \ref{class}.  
         \item For the (a,b)-module \ $E$ \  associated to an holomorphic germ at the origine of \ $\mathbb{C}^{n+1}$ \ with an isolated singularity we have the uniform bounds \ $or(E) \leq n$ \ and \ $\ L(E) \leq n $ \ so the choice \ $N(E) = 2n+ \mu +1$ \ is always possible, where \ $\mu$ \ is the Milnor number (equal to the rank). \ $\hfill \square$
     \end{enumerate}
     
     \parag{Proof} We shall make an induction on the rank of \ $E$. So we shall assume that the result is proved for a rank \ $p - 1$ \ (a,b)-module and we shall consider a regular (a,b)-module \ $E$ \ of rank \ $p \geq 1$, an (a,b)-module \ $E'$, an integer \ $N \geq N_0(E)$ \ and an \ $\A-$isomorphism \ $\varphi$ \ as in \ $(1)$.\\
     From \ref{det. finie ordre reg.} we know that \ $E'$ \ is then regular and has order of regularity \ $or(E') = or(E)$.\\
         Choose now a complex number \ $\lambda$ \ which is minimal \ $modulo \ \mathbb{Z} $ \ such there exists an exact sequence of (a,b)-module ( so \ $\lambda = \tilde{\lambda}_{min}(E)$ \ with the terminology of \S 1.3)
     \begin{equation*}
      0 \to E_{\lambda} \overset{\alpha}{\to}  E \overset{\beta}{\to} F \to 0 . \tag{2}
      \end{equation*}
     This exists from theorem \ref{J-H}. The (a,b)-module \ $F$ \ has rank \ $p-1$ \ and from \ref{3} \ we have \ $N_0(E) \geq N(F^*, -\lambda)$. So we know from \ref{4} that the extension \ $(2)$ \ is determined by the extension
       \begin{equation*}
      0 \to E_{\lambda}/b^N.E_{\lambda} \overset{\alpha_N}{\to}  E/b^N.E \overset{\beta_N}{\to} F/b^N.F \to 0 . \tag{$2_N$}
      \end{equation*}
      Now, using the \ $\A-$isomorphism \ $\varphi$ \ we obain an injective \ $\A-$linear map
      $$ j_N :  E_{\lambda}/b^N.E_{\lambda} \hookrightarrow E'/b^N.E' .$$
      Using the proposition \ref{sub min} with the fact that \ $N \geq or(E') + 2$ \ there exists a unique normal inclusion \ $j : E_{\lambda} \hookrightarrow E' $ \ inducing \ $j_N$.\\
      Define \ $F' : = E'/j(E_{\lambda})$. Then \ $F'$ \ is a rank \ $p-1$ \ (a,b)-module and the exact sequence
      \begin{equation*}
      0 \to E_{\lambda} \overset{j}{\to} E' \to F' \to 0 \tag{2'}
      \end{equation*}
      induced the extension \ $(2_N)$. Using the induction hypothesis, because the inequalities \ $or(E) \geq or(F)$ \ from \ref{det. finie ordre reg.}  and \ $L(E) + rank(E) \geq L(F) + rank(F)$ \ from \ref{maj. lambda(max)(F)} implies \ $N_0(E) \geq N_0(F)$ , we have a unique isomorphism  \ $\Psi : F \to F' $ \ compatible with the one induced by \ $\varphi$ \ between \ $F/b^N.F$ \ and \ $F'/b^N.F'$. Using \ref{4},  \ref{3} and  the inequality \ $N_0(E) \geq N(F^*, -\lambda)$ \  we have an unique isomorphism of extensions
      $$\xymatrix{ 0 \ar[r] & E_{\lambda} \ar[d]^= \ar[r]^{\alpha} & E \ar[d]^{\Phi} \ar[r]^{\beta} & F \ar[d]^{\Psi} \ar[r] & 0 \\
      0 \ar[r] &  E_{\lambda} \ar[r]^j & E' \ar[r] & F' \ar[r] & 0 }$$
      concluding the proof. $\hfill \blacksquare$
      
      \bigskip
      
         \begin{lemma}\label{final} 
         Let \ $E : = J_k(\lambda) $ \ the rank \ $k$ \ (a,b)-module defined by the \ $\mathbb{C}[[b]]-$basis \ $e_1, \cdots, e_k$ \ and by the following relations
      $$ a.e_j = (\lambda + j - 1).b.e_j + e_{j+1} \quad \forall  j \in [1,k]  $$
      with the convention \ $e_{k+1} = 0 $. We have \ $\delta(E) = or(E) = k -1, L(E) = -k+1$. The integer \ $or(E) + L(E) + rank(E) + 1 = k+1$ \ is the best possible for the theorem.
      \end{lemma}
      
      \parag{Proof} It is easy to see that the saturation \ $E^{\sharp}$ \ is generated by \ $e_1, b^{-1}.e_2, \cdots, b^{-k+1}.e_k$. This gives the equality \ $\delta(E) = or(E) = k -1$.\\
      Assume that we have an inclusion \ $E_{\mu} \hookrightarrow E$ \ such that \ $e_{\mu} \not\in b.E$. Then there exists \ $(\alpha_1, \cdots, \alpha_k) \in \mathbb{C}^k\setminus\{0\} $ \ such that
      $$ a.(\sum_{h =1}^k \ \alpha_h.e_h) = \mu.b(\sum_{h=1}^k \alpha_h.e_h) + b^2.E .$$
      Then we obtain
      $$ \sum_{h=1}^k \ \alpha_h.\big((\lambda + j-1).b.e_h + e_{h+1}\big) =  \sum_{h=1}^k \ \alpha_h.\mu.b.e_h + b^2.E $$
      and so \ $\alpha_1=  \cdots = \alpha_{k-1} = 0 $ \ and we conclude that \ $\mu = \lambda + k - 1$.\\
      An easy computation shows that \ $J_k(\lambda)^* = J_k(-\lambda-2k +2)$ \ and so we have \ $\lambda_{max} = \lambda$. So \ $L(E) = -k + 1$.\\
      Now we shall prove that the integer \ $k + 1$ \ is the best possible in the theorem \ref{finite det. thm} for \ $E = J_k(\lambda)$ \ by giving a regular (a,b)-module  \ $F$ \ such that \ $F/b^k.F \simeq E/b^k.E$ \ and not isomorphic to \ $E$.\\
      Let consider the rank \ $k$ \ (a,b)-module \ $F$ \ defined by  \ $ \sum_{j=1}^k \ \mathbb{C}[[b]].e_j$ \ with the following relations
      \begin{align*}
      & a.e_j =  (\lambda + j -1).b.e_j + e_{j+1} \quad \forall j \in [1,k] \\
      & a.e_k = (\lambda + k -1).b.e_k + \sum_{h=1}^{k-1} \ \alpha_h.b^{k-h+1}.e_h 
      \end{align*}
      Then define, for \ $\beta_1, \cdots, \beta_{k-1} \in \mathbb{C}$,
      $$\varepsilon : = e_k + \sum_{j=1}^{k-1} \ \beta_j.b^{k-j}.e_j .$$
      We have
   \begin{align*}
     a.\varepsilon : = &(\lambda + k -1).b.e_k  + \sum_{h=1}^{k-1} \ \alpha_h.b^{k-h+1}.e_h  \ +\\
    \qquad  \qquad   & \sum_{j=1}^{k-1} \beta_j.\big[b^{k-j}.\big((\lambda+j-1).b.e_j + e_{j+1}) + (k-j).b^{k-j+1}.e_j \big] \\
    a.\varepsilon : = &(\lambda + k -1).b.e_k  + \sum_{h=1}^{k-1} \big(\alpha_h + \beta_h.(\lambda + k-1) + \beta_{h-1}\big).b^{k-h+1}.e_h 
   \end{align*}
   
   Let now choose \ $\beta_1, \cdots, \beta_{k-1}$ \ such that we have 
   $$ \alpha_h + \beta_h.(\lambda + k - 1) + \beta_{h-1} = (\lambda + k -1 + \beta_{k-1}).\beta_h \quad \forall h \in [1,k-1]$$
   with the convention \ $\beta_0 = 0 $. We obtain the system of equations
   $$ \alpha_h + \beta_{h-1} = \beta_{k-1}.\beta_h \quad \forall h \in [1,k-1].$$
   This implies, assuming \ $\beta_{k-1} \not= 0$, that \ $\beta_{k-1}$ \ is solution of the equation
   $$ x^k = \alpha_{k-1}.x^{k-2} + \cdots + \alpha_2.x + \alpha_1.$$
   Now choose \ $\alpha_2 = \cdots = \alpha_{k-1} = 0 $ \ and \ $\alpha_1 : = \rho^k $ \ with \ $\rho \in ]0,1[$. Then  choose \ $\beta_j = \rho^{k-j} \quad \forall j \in [1,k-1]$. It is clear that the corresponding \ $F_{\rho}$ \ satisfies \ $F/b^k.F \simeq E/b^k.E$ \ as \ $a.e_k  = e_k + \rho.b^k.e_1$ \ in \ $F_{\rho}$. But the relation \ $a.\varepsilon = (\lambda + k -1 + \rho^k).b.\varepsilon $ \ with \ $\varepsilon \not= 0$ \  shows that \ $F_{\rho}$ \ cannot be isomorphic to \ $J_k(\lambda)$. $\hfill \blacksquare$

      \section*{Bibliography}
      
      \bigskip
      
      \begin{enumerate}
      
      \item{[B.93]} Barlet, Daniel \textit{Theory of (a,b)-modules I} \ in Complex Analysis and Geometry, Plenum Press New York (1993), p.1-43.
      
      \item{[B.95]} Barlet, Daniel \textit{Theorie des (a,b)-modules II. Extensions} in Complex Analysis and Geometry, Pitman Research Notes in Mathematics Series  366 Longman (1997), p. 19-59.
      
      \item{[D.70]} Deligne, Pierre \textit{\'Equations diff\'erentielles \`a points singuliers reguliers} Lect. Notes in Maths, vol. 163, Springer-Verlag (1970).
      
      \end{enumerate}

 \end{document}